%% file: paper_ver9_pgm.tex
\documentclass{aims}
\usepackage{amsmath}
  \usepackage{paralist}
  \usepackage{graphics} %% add this and next lines if pictures should be in esp format
  \usepackage{epsfig} %For pictures: screened artwork should be set up with an 85 or 100 line screen
\usepackage{graphicx}  \usepackage{epstopdf}%This is to transfer .eps figure to .pdf figure; please compile your paper using PDFLeTex or PDFTeXify.
 \usepackage[colorlinks=true]{hyperref}
\hypersetup{urlcolor=blue, citecolor=red}

\def\currentvolume{X}
 \def\currentissue{X}
  \def\currentyear{200X}
   \def\currentmonth{XX}
    \def\ppages{X--XX}

\newtheorem{theorem}{Theorem}[section]

\newtheorem{lemma}[theorem]{Lemma}
\newtheorem{proposition}{Proposition}

\theoremstyle{definition}

\newtheorem{remark}{Remark}

\usepackage{float}
\input{_definitions}

\newcounter{rmnum}
\newenvironment{romannum}{\begin{list}{{\upshape (\roman{rmnum})}}{\usecounter{rmnum}
\setlength{\leftmargin}{22pt}
\setlength{\rightmargin}{8pt}
\setlength{\itemsep}{2pt}
\setlength{\itemindent}{-1pt}
}}{\end{list}}

\newcounter{anum}

\makeatletter
\DeclareRobustCommand*\cal{\@fontswitch\relax\mathcal}
\makeatother

\def\clA{{\cal A}}
\def\clX{{\cal X}}
\def\clY{{\cal Y}}
\def\clZ{{\cal Z}}
\def\E{{\sf E}}
\def\bS{\mathbb{S}}
\def\sP{{\sf P}}
\def\tsP{\tilde{\sf P}}
\def\bpi{{\bar{\pi}}}
\def\Zq{{\hat{V}}}
\def\barh{{\bar{h}}}
\def\barg{{\bar{g}}}

\title[FPF for Collective Inference] %Use the shortened version of the full title
      {Feedback Particle Filter for Collective Inference}

\author[Jin-Won Kim, Amirhossein Taghvaei, Yongxin Chen, and Prashant G. Mehta]{}

\subjclass{Primary: 60G35, 62M20; Secondary: 94A12.}
 \keywords{Collective inference, Feedback particle filter.}

 \email{jkim684@illinois.edu}
 \email{amirhoseintghv@gmail.com}
 \email{yongchen@gatech.edu}
 \email{mehtapg@illinois.edu}

\thanks{Kim and Mehta are supported in part by the C3.ai Digital
          Transformation Institute sponsored by C3.ai Inc. and the
          Microsoft Corporation, and in part by the National Science
          Foundation grant NSF 1761622.  Chen is supported by
          the NSF 2008513.}

\thanks{$^*$ Corresponding author: Prashant G. Mehta}

\begin{document}
\maketitle

% Enter the first author's name and address:
\centerline{\scshape Jin-Won Kim}
\medskip
{\footnotesize
% please put the address of the first author
 \centerline{Coordinated Science Laboratory}
   \centerline{University of Illinois, Urbana-Champaign }
   \centerline{1308 W. Main St.}
   \centerline{ Urbana, IL 61801, USA}
} % Do not forget to end the {\footnotesize by the sign }

\medskip

\centerline{\scshape Amirhossein Taghvaei}
\medskip
{\footnotesize
 % please put the address of the second  and third author
 \centerline{Department of Mechanical and Aerospace Engineering }
   \centerline{University of California, Irvine }
   \centerline{4100 Calit2 Building}
   \centerline{Irvine, CA 92697-2800, USA}
}
\medskip

\medskip

\centerline{\scshape Yongxin Chen}
\medskip
{\footnotesize
 % please put the address of the second  and third author
 \centerline{School of Aerospace Engineering}
   \centerline{Georgia Institute of Technology}
   \centerline{Guggenheim 448B}
   \centerline{Atlanta, GA 30332, USA}
}
\medskip

\centerline{\scshape Prashant G. Mehta$^*$}
\medskip
{\footnotesize
 % please put the address of the second  and third author
 \centerline{Coordinated Science Laboratory}
   \centerline{University of Illinois, Urbana-Champaign }
   \centerline{1308 W. Main St.}
   \centerline{Urbana, IL 61801, USA}
}

\bigskip

% The name of the associate editor will be entered by an editorial staff
% "Communicated by the associate editor name" is not needed for special issue.
 \centerline{(Communicated by the associate editor name)}

%The abstract of your paper
\begin{abstract}
The purpose of this paper is to describe the feedback particle filter
algorithm for problems where there are a large number ($M$) of
non-interacting agents (targets) with a large number ($M$) of non-agent
specific observations (measurements) 
that originate from these agents.  In its basic form, the problem is
characterized by data association uncertainty whereby the association
between the observations and agents must be deduced in addition to 
the agent state.  In this paper, the large-$M$ limit is interpreted as
a problem of collective inference.  This viewpoint is used to derive the
equation for the empirical distribution of the hidden agent states.  A
feedback particle filter (FPF) algorithm for this problem is presented and
illustrated via numerical simulations.  Results are presented for the
Euclidean and the finite state-space cases, both in continuous-time
settings.  The classical FPF algorithm is shown to be the special case
(with $M=1$) 
of these more general results.  The simulations help show that
the algorithm well approximates
the empirical distribution of the hidden states for large $M$.   
\end{abstract}

%%%%%%%%%%%%%%%%%%%%%%%%%%%%%%%%%%%%%%%%%%%%%%%%%%%%%%%%%%%%%%%%%%%%%%%%%%%%%%

\section{Introduction}
\label{sec:intro}

Filtering with data association uncertainty is important to a number
of classical applications, including target tracking, weather surveillance,
remote sensing, autonomous navigation and
robotics~\cite{bar2009probabilistic,thrun2002probabilistic}.    
Consider,
e.g., the problem of multiple target tracking with radar. The
targets can be multiple aircrafts in air defense, or multiple weather
cells in weather surveillance, or multiple landmarks in autonomous
navigation and robotics.  In each of these applications, there exists data
association uncertainty in the sense that one can not assign, in an
apriori manner, individual observations (measurements) to individual
targets.  
Given the large number of applications, algorithms for filtering
problems with data association uncertainty have been extensively
studied in the past; cf.,~\cite{bar2009probabilistic,kirubarajan2004probabilistic,kyriakides2008sequential} %,survey_data_assoc}
and references therein.  The feedback particle filter (FPF) algorithm for
this problem appears in~\cite{yang2018probabilistic}.  

The filtering problem with data association uncertainty is 
closely related to the filtering problem with aggregate and anonymized
data.  Some of these problems have gained in importance recently
because of COVID-19.   Indeed, the
spread of COVID-19 involves dynamically evolving hidden processes (e.g.,
number of infected, number of asymptomatic etc..) that must be
deduced from noisy and partially observed data (e.g., number of tested
positive, number of deaths, number of hospitalized etc.).  
In carrying out data assimilation for such problems, one typically only has
aggregate observations.  For
example, while the number of daily tested positives is available, the
information on the disease status of any particular agent in the
population is not known.  Such problems are referred to as collective
or aggregate inference problems~\cite{sheldon2013approximate,sheldon2011collective,chen2018state}.

In a recent important work~\cite{singh2020inference}, algorithms are
described for solving the collective inference problem in 
graphical models, based on the large deviation theory.  These results
are also specialized to the smoothing problems for the hidden Markov
models (HMMs).  Two significant features of these algorithms are: 
(i) the complexity of the data assimilation does not grow with
the size of the population; and (ii) for a single agent, the algorithm
reduces to the classical forward-backward smoothing algorithm for
HMMs.

The main purpose in this paper is to interpret the collective inference
problem as a limit of the data association problem, as the number of
agents ($M$) become large.  Indeed, for a small number of agents, data
association can help reduce the uncertainty and improve the
performance of the filter.  
However, as number of agents gets larger, the data association-based
solutions become less practical and may not offer much benefit: On the
one hand, the complexity of the filter grows because the number of
associations for $M$ agents with $M$ observations is $M!$.  On the
other hand, the performance of any practical algorithm is expected to
be limited.  

In this paper, the filtering problem for a large number of agents is
formulated and solved as a collective inference problem.  Our main
goal is to develop the FPF algorithm to solve the
collective filtering problem in continuous-time settings.
For this purpose, the Bayes' formula for collective inference is introduced
(following the formulation of~\cite{singh2020inference}) and
compared with the standard Bayes' formula.  The collective Bayes'
formula is specialized to derive the equations of collective filtering
in continuous-time settings.  The FPF algorithm for
this problem is presented for the
Euclidean and the finite state-space cases.  The classical FPF and
ensemble Kalman filter (see~\cite{bergemann2012, surace2019avoid,
  reich2011dynamical,reich2015probabilistic}) algorithms are shown to
be the special case (with $M=1$) of these more general results.  The
algorithm is illustrated using numerical simulations.

The outline of the remainder of this paper is as follows.  The problem
formulation appears in Section~\ref{sec:prelim}. The collective Bayes
formula and the filter equations are derived in
Section~\ref{sec:collective-bayes-formula}.  The FPF algorithm is
described in Section~\ref{sec:CF-FPF}.  The simulation results appear
in Section~\ref{sec:sims}.  All the proofs are contained in the Appendix.

\section{Problem formulation}\label{sec:prelim}

%TASKS FOR JINKIM:
%change the language from "target" to "agents"
%make the time axis more general - like you do in your version
%Include examples for the linear and HMM.

\subsection{Dynamic model}

\begin{enumerate}
\item There are $M$ agents. The set of agents is denoted as ${\cal M}
  = \{1,2,\ldots, M\}$.  The set of permutations of ${\cal M}$ is
  denoted by $\Pi({\cal M})$ whose cardinality $|\Pi({\cal M})| =
  M!$.

\smallskip

\item Each agent has its own state--observation pair. The state
  process for the $j^{\text{th}}$ agent is $X^j = \{X_t^j:t\in I\}$, a
  Markov process and $I$ is the index set (time). The associated
  observation for the $j^{\text{th}}$ agent is $Z^j =
  \{Z_t^j:t\in I\}$.%  is modeled for each time $t \in I$ by the observation
  % kernel for ${\sf P}(Z_t\mid X_t)$.

\smallskip

\item At time $t$, the observations from all agents is aggregated
  while their labels are censored through random permutations $\sigma_t
  \in \Pi({\cal M})$. The association $\sigma_t =
  (\sigma_t^1,\sigma_t^2,\ldots,\sigma_t^M)$ signifies that the
  $j^{\text{th}}$ observation originates from agent $\sigma_t^j$.  The
  random permutations is modeled as a Markov process on $\Pi({\cal
    M})$.  

% \item The centralized estimator is provided with the aggregated observation without the true association $\sigma_t$. Therefore, it does not have access on individual observation $Z^j$, but only the empirical distribution is available.

%\item The estimator seeks to estimate the distribution of $\{X_t^j:j\in {\cal M}\}$ given the aggregated observations. 
\end{enumerate}

\subsection{Standard filtering problem with data association}

In its most general form, the filtering problem is to assimilate the
measurements $\clZ_t=\sigma(Z_t^i:1\leq i \leq M, 0\leq s \leq t)$ to
deduce the posterior distribution of the hidden states $\{X_t^i:1\leq
i \leq M\}$.  Given the associations are also hidden, the problem is
solved through  building a filter also to estimate the permutation
$\sigma_t$.

A number of approaches have been considered to solve the problem in a
tractable fashion:  Early approaches included multiple hypothesis
testing (MHT) algorithm, requiring exhaustive
enumeration~\cite{reid1979algorithm}. However,
exhaustive enumeration leads to an NP-hard problem because
number of associations increases exponentially with time.  The
complexity issue led to development of the probabilistic MHT or its
simpler ``single-scan'' version, the joint probabilistic data
association (JPDA) filter~\cite{kirubarajan2004probabilistic, Bar-ShalomYaakov1988Tada}.  These algorithms are based on computation (or
approximation) of the {\em observation-to-target association
  probability}.  The feedback particle filter extension of the JPDA
filter appears in~\cite{yang2018probabilistic}.

\subsection{Collective filtering problem}

In the limit of large number of non-agent specific observations, it is
more tractable to consider directly the empirical distribution of the
observations: 
\begin{align*}
q_t(z) := \frac{1}{M}\sum_{j=1}^M \delta_{Z^j_t}(z)
\end{align*}
and use it to estimate the empirical distribution of the hidden
states -- denoted as $\pi_t$ at time $t$.  The problem is referred to
as the  {\em collective filtering} problem.  
%
% The problem of {\em collective filtering} is to estimate the
% distribution of entire agents given the distribution of the
% observations from each agent. 
%Following the work of \cite{singh2020inference}, the remainder of this paper will describe the mathematical objective.

%the optimal estimate is the solution to the minimization problem:
%\begin{subequations}\label{eq:KL-optimization}
%\begin{align}
%\mathop{\text{Minimize :}}_{\tsP\in {\cal P}(X,Z)}\quad &{\sf D}\big(\tsP(X,Z)\|\sP(X,Z)\big) \label{eq:KL-optimization-a}\\
%\text{Subject to:} \quad& \tsP|_{Z_s} = q_s, \quad \forall s\in I,\; s\le \tau \label{eq:KL-optimization-b}
%\end{align}
%\end{subequations}
%where ${\cal P}(X,Z)$ is the family of joint probability distribution on both $X$ and $Z$. The constraint is such that the marginal on $Z$ follows the given information $q$. The optimal distribution is then the marginal of $\tsP$ on $X$.
%
%Let $\pi_{t\mid \tau}$ denote the marginal of the optimal $\tsP$ on $X_t$, given information of $q$ up to time $\tau\in I$.

\section{Collective Bayesian filtering}\label{sec:collective-bayes-formula}

It is easiest to introduce the set-up in discrete time finite
state and observation space case.  This is done
prior to presenting the generalization to the continuous-time model.  

\subsection{Discrete time case}    %One-step estimator for HMM
% In this section, we consider discrete time setting. 
The 
index set $I = \{0,1,\ldots\}$. The state space $\clX =
\{1,2,\ldots, d\}$ and the observation space ${\cal Y} =
\{1,2,\ldots, m\}$ are assumed to be both finite.
The system is modeled by the state transition matrix $\sP(X_{t+1}=x\mid X_t =
x') = p(x\mid x')$ and the emission matrix $\sP(Z_t=z\mid X_t=x) =
o(z\mid x)$.  The set of joint probability distributions on $\clX \times \clY$ is
denoted as ${\cal P}(X,Z)$.

% In the following, $\pi_t$ is the empirical distribution of $X_t$
% and $\pi_{t+1|t}(x^+) := \sum_{x} p(x^+\mid x) \pi_t(x) $.  

\newP{Optimization problem} Given $\pi_t$ and $q_{t+1}$ the one-step collective
inference problem is 
% Let $\pi_{t\mid \tau}$ be the
% (estimated) empirical distribution of $X_t$ given the empirical
% distribution of the observation $q_{\cdot}(z)$ up to
% time $\tau$. Given $\pi_{t\mid t}$, consider the following optimization problem to obtain the optimal one-step estimate:
\begin{subequations}\label{eq:KL-optimization}
\begin{align}
\mathop{\text{Minimize :}}_{\tsP\in {\cal P}(X,Z)}\quad &{\sf
                                                          D}\big(\tsP \mid
                                                          \sP \big) \label{eq:KL-optimization-a}\\
\text{Subject to:} \quad& \sum_{x} \tsP (x,z)  = q_{t+1}(z) \label{eq:KL-optimization-b}
\end{align}
\end{subequations}
where $\sP(x,z) =  \sum_{x'} p(x\mid x') \pi_t(x')   \, o(z|x)$ and ${\sf D}(\cdot|\cdot)$ is the K-L divergence.
%where $\sP(x,z) =  \pi_{t+1|t}(x)  \, o(z|x)$.  
% where ${\cal P}(X,Z)$ is the family of joint probabilities on $\clX
% \times \clY$, and $\sP(X_{t+1},Z_{t+1})$ is the nominal estimate on
% $X_{t+1}$ and $Z_{t+1}$ given $\pi_{t\mid t}$.

\medskip

The justification for considering this type of optimization objective
to model the collective inference problem is given
in~\cite[Sec. III-A]{singh2020inference}.  The K-L divergence is the rate
function in the large deviation theory that characterizes the decay in
probability of observing an empirical distribution~\cite{chen2016relation}.  The
solution to this problem is described in the following proposition
whose proof appears in Appendix~\ref{apdx:one-step-estimator}.

\begin{proposition}\label{prop:one-step}
Consider the optimization problem~\eqref{eq:KL-optimization}. The optimal one-step estimate is given by:
\begin{subequations}\label{eq:one-step-estimator}
\begin{align}
\pi_{t+1\mid t}(x) &= \sum_{x'\in\clX} p(x\mid x')\pi_{t\mid t}(x') \label{eq:one-step-estimator-a}\\
\pi_{t+1\mid t+1}(x) &= \sum_{z\in{\cal Y}} \frac{o(z\mid x)\pi_{t+1\mid t}(x)}{\xi_{t+1}(z)}q_{t+1}(z) \label{eq:one-step-estimator-b}\\
\xi_{t+1}(z) &= \sum_{x\in\clX} o(z\mid x)\pi_{t+1\mid t}(x) \label{eq:one-step-estimator-c}
\end{align}
\end{subequations}
\end{proposition}

%\begin{subequations}\label{eq:one-step-estimator}
%\begin{align}
%\pi_{t+1\mid t}(x_{t+1}) &= \int p(x_{t+1}\mid x_t)\pi_{t\mid t}(\ud x_t) \label{eq:one-step-estimator-a}\\
%\pi_{t+1\mid t+1}(x_{t+1}) &\propto \int \frac{p(z_{t+1}\mid x_{t+1})\pi_{t+1\mid t}(x_{t+1})}{\xi_{t+1}(z_{t+1})}q_{t+1}(\ud z_{t+1}) \label{eq:one-step-estimator-b}\\
%\xi_{t}(z_t) &= \int p(z_t\mid x_{t})\pi_{t+1\mid t}(\ud x_{t}) \label{eq:one-step-estimator-c}
%\end{align}
%\end{subequations}

%\medskip
%
%If state and observation take value in finite space, it is simplified to the following (cf. Yongxin's algorithm for smoothing problem where $\tau = T$).
%\begin{align*}
%\pi_{t+1\mid t}(x_{t+1}) &= \sum_{x_{t}} p(x_{t+1}\mid x_t)\alpha_t(x_t)\gamma_t(x_t)\\
%\gamma_t(x_t) &= \sum_{z_t} p(z_t\mid x_t) \frac{q_t(z_t)}{\xi_t(z_t)}\\
%\xi_t(z_t) &= \sum_{x_t} p(z_t\mid x_t)\pi_{t\mid t-1}(x_t)
%\end{align*}
%

\begin{remark}
With $M=1$, the one-step estimate~\eqref{eq:one-step-estimator}
reduces to the Bayes' formula.
\end{remark}

\begin{remark}
The optimization problem~\eqref{eq:KL-optimization} is a special case
of the problem introduced in~\cite{singh2020inference} over a time
horizon. 
\end{remark}

\subsection{Continuous time case} \label{sec:ctmc}

The recursive formula for the discrete time problem is extended to the
continuous-time problem. In this case, the index set $I = [0,\infty)$
and the state process $\{X_t:t\ge 0\}$ is a continuous-time Markov
process with the generator denoted as $\clA$.  The associated adjoint
is denoted as $\clA^\dagger$.  Additional details on the particular
types of Markov processes appears in the sequel.  
% Three special cases of interest are:
% (i) Ito diffusions in the Euclidean state-space, (ii) the linear
% Gaussian case, and (iii) the finite state-space case.  
% In the simplest settings, $\clX$ is finite and
% $\clA$ is the probability transition rate $r(x\mid x')$ from state
% $x'$ to $x$; in this case $\clA^\dagger$ is its matrix transpose. 

We assume the following stochastic differential equation (SDE) model
for scalar-valued observations:
\begin{equation}\label{eq:obs_model}
\ud Z_t = h(X_t) \ud t + \sigma_w \ud W_t
\end{equation}
where $\{W_t:t\ge 0\}$ is the standard Wiener process and $\sigma_w>0$. It is assumed
that the observation noise is independent of the state process.

In the continuous-time settings, the empirical distribution
$q=\{q_t:t\geq 0\}$ of the
observations is defined for an increment $\Delta Z_t^j := Z_{t+\Delta t}^j -
Z_t^j$ for $j=1,2,\hdots,M$. We denote the mean process by:
\begin{equation*}
\Delta \hat{Z}_t := \int \Delta z \,\ud q_t(\Delta z) = \frac{1}{M}
\sum_{j=1}^M \Delta Z_t^j 
\end{equation*}
The following quantities related to the second-moment are also of interest:
\begin{align*}
\Zq_t := & \lim_{\Delta t \to 0} \frac{1}{\Delta t}\big(\Delta \hat{Z}_t\big)^2\\
V_t :=& \lim_{\Delta t \to 0} \frac{1}{\Delta t} \int (\Delta z)^2
        \,\ud q_t(\Delta z) - \Zq_t = \lim_{\Delta t \to 0}
        \frac{1}{M\Delta t}
\sum_{j=1}^M (\Delta Z_t^j)^2  - \Zq_t
\end{align*}
%
%first and the second moment by:
%\begin{align*}
%\int \Delta z \,\ud q_t(\Delta z) &=: \Delta \hat{Z}_t\\
%\int (\Delta z)^2  \,\ud q_t(\Delta z) & =: \Delta [Z]_t
%\end{align*}
%The interpretation for these definitions is as follows: $\hat{Z}_t$ is the
%mean of the observation processes. As $\Delta t \to 0$, the second moment becomes increment of the quadratic variation $[Z]_t$. Assume that the quadratic variation exists and also it is differentiable. Define
%\begin{align*}
%\Zq_t :=& \frac{\ud}{\ud t}[\hat{Z}]_t\\
%V_t :=& \frac{\ud}{\ud t} \big([Z]_t - [\hat{Z}]_t\big)
%\end{align*}
It is assumed that these limits are well-defined.  The first limit
represents the quadratic variation of the stochastic process
$\{\hat{Z}_t:t\geq 0\}$ and the second limit is an empirical variance
of the increments. 
% These terms represent two factors: The quadratic variation of $\hat{Z}_t$ itself; and the variance of the observations around its mean, respectively.

\newP{Optimal recursive update} Using this notation, the continuous
time counterpart of the Prop.~\ref{prop:one-step} is as follows.  The
proof appears in the Appendix~\ref{apdx:continuous-time}.

\medskip

\begin{proposition}\label{prop:continuous-time}
%For each $x\in \clX$, the optimal estimate is the following:
The collective filter $\{\pi_t:t\geq 0\}$ solves the following
evolution equation:
\begin{align}
\ud \pi_t(x) = {\cal A}^\dagger \pi_t (x)\ud t  & + \frac{1}{\sigma_w^2}\pi_t(x)(h(x)-\hat{h}_t)(\ud \hat{Z}_t - \hat{h}_t \ud t) \label{eq:cont-time-update}\\
&+ \frac{1}{\sigma_w^2}
\pi_t(x)\big(g(x) - \hat{g}_t \big)\Big(\frac{V_t+\Zq_t}{\sigma_w^2}-1 \Big)\ud t \nonumber
\end{align}
where $\hat{h}_t = \int h(x)\pi_t(x)\ud x$, $g(x)=\frac{1}{2}(h(x)-\hat{h}_t)^2$ and $\hat{g}_t = \int g(x) \pi_t(x)\ud x$. 
\end{proposition}

\medskip

\begin{remark}
The  first two terms of the righthand-side
of~\eqref{eq:cont-time-update} are identical to the Wonham filter.
The third term is an additional correction term that arises on account
of the second moment of $q$. A special case ($M=1$) is when the empirical
distribution $q_t$ is a Dirac delta measure concentrated on a single trajectory
$\hat{Z}_t=Z_t$.  In this case the variance $V_t = 0$ and the third
term vanishes (using the Ito's rule $\ud \hat{Z}_t^2 = \Zq_t\ud t =
\sigma_w^2\ud t$).  Therefore, for $M=1$, the collective filter
reduces to the Wonham filter.  

Another special case is when $M=\infty$.  In
this case, $\hat{Z}_t$ is a deterministic process $\int_0^t
\E\big(h(X_s)\big)\ud s$ and $V_t = \sigma_w^2$ almost surely. In this
case, both the second and third terms
of~\eqref{eq:cont-time-update} are zero and the collective filter reduces
to the Kolomogorov's forward equation.
\end{remark}

\subsection{Linear-Gaussian case}

A special case of the continuous-time model is the linear-Gaussian case where the state and observation processes are defined by:
\begin{subequations}\label{eq:CLG-system}
\begin{align}
\ud X_t &= AX_t \ud t + Q^{1/2}\ud B_t,\quad X_0 \sim {\cal N}(m_0,\Sigma_0)\\
\ud Z_t &= HX_t \ud t + \sigma_w \ud W_t
\end{align}
\end{subequations}
Here the drift term and the observation function are linear and the noise processes and initial condition are Gaussian.%  In this case, FPF approximates Kalman-Bucy filter~\cite{taghvaei2020optimal}. For continuous-time problems involving data association
% uncertainty, a Markov model of $\{\sigma_t\}_{t\geq 0}$ appears in~\cite{yang2018probabilistic}.
% Without data association, the observations from multiple agents is mixture of $M$ Gaussians. In order to compute obtain Gaussian estiate, $q_t$ is approximated by a Gaussian with mean $\Delta \hat{Z}_t$ and variance $V_t \Delta t$.  

%where $\{B_t\}_{t\geq 0}$, $\{W_t\}_{t\geq 0}$ are standard Wiener
%processes, $Q$ and $R$ are noise variances. It is assumed
%that $X_0$, $\{B_t\}_{t\geq 0}$, and $\{W_t\}_{t\geq 0}$ are mutually
%independent.  For continuous-time problems involving data association
%uncertainty, a Markov model of $\{\sigma_t\}_{t\geq 0}$ appears in~\cite{yang2018probabilistic}.

%%where $\hat{z}_t = \frac{1}{M}\sum_{j=1}^M Z_t^j$ is the empirical
%%mean and $\tilde{W}_t^j$ is a standard Wiener process introduced to
%%`fit' the empirical distribution to the observations.  The coefficient
%%$V_t$ has the meaning of the quadratic variation of the random
%%component. 
%
%% system of interacting particles to approximate Kushner equation by its empirical distribution~\cite{yang2016}. For linear-Gaussian problem, it approximates Kalman-Bucy filter~\cite{taghvaei2020optimal}. The goal of this section is to derive collective filtering equation for continuous time linear Gaussian problem, and then obtain linear FPF that approximate the collective filter.

The following is the counterpart of~\eqref{eq:cont-time-update} for the linear-Gaussian case.  Its proof appears in the Appendix~\ref{apdx:pf-continuous-LG}.  

\medskip

\begin{proposition}\label{prop:continuous-LG} 
Consider the collective filtering problem for the continuous-time linear
Gaussian model~\eqref{eq:CLG-system}.  Then $\pi_{t}$ is a Gaussian whose mean $m_{t}$ and variance $\Sigma_{t}$ evolve according to
\begin{subequations}\label{eq:ckf}
\begin{align}
\ud m_{t} &= Am_{t} \ud t + \frac{1}{\sigma_w^2}\Sigma_{t}H^\top \big(\ud \hat{Z}_t - Hm_{t} \ud t\big)\\
\frac{\ud}{\ud t} \Sigma_{t} &= A\Sigma_{t} + \Sigma_{t} A^\top + Q - \frac{1}{\sigma_w^2}\Sigma_{t} H^\top \big(1 - \frac{V_t}{\sigma_w^2}\big) H \Sigma_{t}
\end{align}
\end{subequations}
with the initial conditions specified by $m_0$ and $\Sigma_0$.
\end{proposition}

\section{Collective Feedback Particle Filter}\label{sec:CF-FPF}

%In this section,  FPF algorithm is presented in order to obtain numerical solution to~\eqref{eq:cont-time-update}. 
The feedback particle filter is a controlled
interacting particle system to approximate the solution of~\eqref{eq:cont-time-update}. 
% the nonlinear filtering problem in continuous-time
% settings~\cite{yang2016}.
% The construction of the FPF is in two
% steps:  
% In step 1, a mean-field system is designed such that the distribution
% of the particles achieves the optimal estimate obtained by~\eqref{eq:cont-time-update}. Secondly, the mean-field
% system is approximated by finite number of particles and their empirical distribution.
In this section, we extend the FPF formulae in~\cite{yang2016, yang2016MarkovChain} to the collective filtering settings.

\subsection{Euclidean setting}

% The original form of the FPF is developed for It\^{o} diffusion
% processes on the state space $\clX=\Re^d$. 
The
 state process is defined by the following SDE:
 \begin{equation*}\label{eq:dyn_sde}
 \ud X_t = a (X_t) \ud t + \sigma(X_t) \ud B_t, \quad X_0\sim \pi_0
 \end{equation*}
 where $a\in C^1(\Re^d; \Re^d)$, $\sigma\in C^2(\Re^d; \Re^{d\times p})$ and $B=\{B_t:t\ge 0\}$ is a standard Wiener process.  $\pi_0$ denotes the initial
 distribution of $X_0$. It is assumed that $\pi_0$ has a probability density with respect to the Lebesgue measure. For~\eqref{eq:obs_model}, the observation function $h\in C^2(\Re^d;\Re)$. % Let ${\cal A}$ denote the infinitesimal generator of the state process, and ${\cal A}^\dagger$ be the adjoint operator.

%For the Euclidean setting, $\pi_t$ is now probability density, and the update for the filter is direct extension of~\eqref{eq:cont-time-update}.
% \begin{align}
% \ud \pi_t(x) = & {\cal A}^\dagger \pi_t (x) + \frac{1}{\sigma_w^2}\pi_t(x)(h(x)-\hat{h}_t)(\ud\hat{Z}_t - \hat{h}_t \ud t) \label{eq:Euc-update} \\+&\frac{1}{\sigma_w^2}
% \pi_t(x)(g(x) - \hat{g}_t  - h(x)\hat{h}_t + \hat{h}_t^2)(\frac{V_t+\Zq_t}{\sigma_w^2}- 1)\ud t \nonumber
% \end{align}

\newP{Feedback particle filter} The formula for collective FPF is as follows:
\begin{align}
\ud X_t^i  =  a(X_t^i)\ud t + \sigma(X_t^i) \ud B_t^i &+ K_t(X_t^i) \circ \big(\ud \hat{Z}_t - \alpha_t h(X_t^i) \ud t - (1-\alpha_t) \barh_t \ud t \Big) \nonumber\\
 &+ \Big(\frac{V_t+\Zq_t}{\sigma_w^2}- 1\big)u_t(X_t^i) \ud t, \quad X_0^i \sim \bpi_0 \label{eq:FPF-euc}
\end{align}
where $B^i = \{B_t^i:t\ge 0\}$ is an independent copy of the process
noise and $\alpha_t := \frac{1}{2}(1- \frac{V_t}{\sigma_w^2})$. The
symbol $\circ$ means that the SDE is expressed in its Stratonovich form. % $\barh_t$ denotes the expectation of $h$ over particles $X_t^i$.
The gain function $K_t(\cdot)$ and $u_t(\cdot)$ are chosen to solve Poisson equations:
\begin{subequations}\label{eq:Poisson-eqns}
\begin{align}
- \nabla \cdot(\bpi_t K_t) &= \frac{1}{\sigma_w^2}\bpi_t (h-\barh_t)\\
- \nabla \cdot(\bpi_t  u_t)  &=\frac{1}{2} \nabla \cdot(\bpi_t K_t (h-\barh_t))  + \frac{1}{\sigma_w^2}\bpi_t(g-\barg_t)
\end{align}
\end{subequations}
where $\bpi_t$ denotes the distribution of $X_t^i$ and $\barh_t=\int
\bpi_t(x) h(x) \ud x$, $\barg_t=\int
\bpi_t(x) g(x) \ud x$.
The following proposition states the exactness property of the
filter. (That is, the distribution of the particles exactly
matches the distribution $\pi_t$ the collective filter). The proof appears in Appendix~\ref{apdx:pf-FPF-Euc}.

\begin{proposition}\label{prop:FPF-Euc}
Consider the FPF~\eqref{eq:FPF-euc}. Then its probability density function $\bpi_t$ evolves according to the equation~\eqref{eq:cont-time-update}. Consequently, if $\bpi_0= \pi_0$, then
\begin{equation*}
\bpi_t = \pi_t \quad \forall t \ge 0
\end{equation*}
\end{proposition}

\medskip

% Note that FPF~\eqref{eq:FPF} is
% not practical because it requires the distribution $\bpi_t$.  In practice, these
% are approximated by finite $N$ number of particles. In this case, for $i = 1,\ldots, N$, $X_0^i$ is $N$ iid samples from $\bpi_0$ and noise processes $B^i$ is also $N$ copies of independent Wiener processes. In this case,
% 	\begin{equation*}
%   	\barh_t^{(N)}:=\frac{1}{N}\sum_{i=1}^N h(X^i_t),\quad
%   	\barg_t^{(N)}
%   	:=\frac{1}{N}\sum_{i=1}^N g(X^i_t)
%   	\end{equation*}
% and the Poisson equations are also accordingly constructed.  In the limit as $N\to\infty$, the empirical distribution converges
% to $\bpi$ because of the propagation of chaos property
% of the 
% interacting particle systems~\cite{taghvaei2020optimal}. 

\subsection{Linear-Gaussian case}

In linear-Gaussian case~\eqref{eq:CLG-system}, the distribution is
Gaussian and completely characterized by its mean $\bar{m}_t$ and the variance
$\bar{\Sigma}_t$. The explicit formulae for the solution of the Poisson
equations appear in the following lemma.
%In its basic form, the feedback particle filter (FPF) is a controlled
%interacting particle system algorithm to approximate the solution of
%the nonlinear filtering problem in continuous-time settings~\cite{yang2016}.
%Although the algorithm was developed for the general nonlinear case,
%here we only consider the linear Gaussian
%model.

\begin{lemma}

In linear-Gaussian case, $K_t(\cdot) = \frac{1}{\sigma_w^2} \bar{\Sigma}_t H^\top$ and $u_t(\cdot) = 0$ solves the Poisson equations~\eqref{eq:Poisson-eqns}. 

\end{lemma}

The proof is straightforward by directly
using~\eqref{eq:Poisson-eqns}.  Consequently, 
the FPF algorithm for the linear-Gaussian case is as follows:
%The feedback particle filter (FPF) algorithm is designed to sample
%from the posterior $\pi_t$.  For the collective filtering problem, FPF
%is comprised of $N$ particles
%$\{X_t^i:1\leq i \leq N, t \ge 0\}$ which are defined according to the
%mean-field model:
\begin{align}
\ud X_t^i &= A X_t^i \ud t + \ud \bar{B}_t^i + \bar K_t\big (\ud \hat{Z}_t-
            (\alpha_tHX_t^i + (I-\alpha_t)H \bar{m}_t)\ud t\big) \nonumber\\
X_0^i &\stackrel{\text{iid}}{\sim} {\cal N}(m_0,\Sigma_0)\label{eq:FPF}
\end{align}
where $\{\bar{B}_t^i: 1\leq i \leq N, t \ge 0\}$
are $N$ independent copies of the process noise $\{B_t\}_{t\geq 0}$, $\bar{K}_t = \frac{1}{\sigma_w^2} \bar{\Sigma}_t H^\top$, $\alpha_t =
\frac{1}{2} (I - \frac{V_t}{\sigma_w^2})$. It is readily seen that the
mean and the variance of the particles evolves exactly according
to~\eqref{eq:ckf}.

%, and $\bar{m}_t$ and $\bar{\Sigma}_t$ are the mean and the variance, respectively, of $X_t^i$.  
%
%The following proposition states the exactness property of the
%filter. (That is, the mean and variance of the particles exactly
%match the mean and the variance of the collective filter).  Its
%proof appears in the Appendix~\ref{apdx:pf-FPF}. 
%
%
%\medskip
%
%\begin{proposition}\label{prop:FPF}
%Consider the FPF~\eqref{eq:FPF}.  Then its mean $\bar{m}_t$ and
%variance $\bar{\Sigma}_t$ evolve
%according to equations for $m_t$ and $\Sigma_t$, respectively.
%Consequently, if $\bar{m}_0=m_0$ and $\bar{\Sigma}_0=\Sigma_0$ then
%\[
%m_t=\bar{m}_t,\quad \bar{\Sigma}_t=\Sigma_t\quad  \forall \;t>0
%\]
%\end{proposition}

\subsection{Finite-state case}  Consider the continuous-time  Markov chain $\{X_t:t\geq 0\}$ defined on the finite state space case $\mathcal X :=\{e_1,e_2,\ldots,e_d\}$ where $e_k$ are standard bases on $\Re^d$ for $k=1,\ldots,d$.  The dynamics of the Markov chain is given by
\begin{equation*}
\ud X_t = \sum_{x \in \mathcal X,x\neq X_t}(x-X_t) \ud \zeta_t^{x,X_t}
\end{equation*}
where $\zeta_t^{x,y}$ is a Poisson processes with rate $r_{x,y}$ for
$x,y \in \mathcal X$. A count from Poisson process $\zeta_t^{x,X_t}$
causes the  Markov-chain to jump to state $x$.    The observation
model is~\eqref{eq:obs_model}.  The FPF update law   for the
finite-state case is as follows:
%%Another noticeable example is a Markov chain on a finite state-space case, as described earlier in~\ref{sec:ctmc}. Assume the Markov chain The FPF algorithm for Markov chain appears in~\cite{yang2016MarkovChain}. The control term is additional jump process according to the observation. (explain?)
%For collective filter~\eqref{eq:cont-time-update}, the particles evolve according to: 
\begin{align}
\ud X^i_t &= \sum_{x\in \mathcal X ,x\neq X^i_t} (x - X^i_t) \ud \zeta^{x,X^i_t}_t + \sum_{x\in \mathcal X,x\neq X^i_t} (x - X^i_t) \ud  \ell^{x,X^i_t}_t \nonumber \\
&+\sum_{x\in \mathcal X ,x\neq X^i_t} (x - X^i_t) \ud  \tilde{\ell}^{x,X^i_t}_t,\quad X^i_0 \sim \pi_0 \label{eq:FPF-finite}
\end{align}
%and $K(\cdot)$, $\tilde{K}(\cdot)$
%such that for any function $f$ 
%\begin{equation}
%\ud f(X^i_t) = \sum_{x\in S} (f(x) - f(X^i_t)) \ud N _{x,X^i_t}
%\end{equation}
where  
%\begin{align*}
%
%\end{align*}
$\ell^{x,y}_t,\tilde{\ell}^{x,y}_t$ are time-modulated Poisson processes of the following form
\begin{align*}
\ell^{x,y}_t = N^{x,y}(U^{x}_t), \quad \tilde{\ell}^{x,y}_t = \tilde{N}^{x,y}(\tilde{U}^{x}_t)
\end{align*}
Here, $N^{x,y}(\cdot)$ and $\tilde{N}^{x,y}(\cdot)$ are standard Poisson processes with rate equal to one.  The inputs $U^x_t$ and $\tilde{U}^x_t$ are defined according to 
\begin{align*}
\ud U^x_t = K_t(x) (\ud \hat{Z}_t- \barh_t \ud t),\quad \ud \tilde{U}^x_t = \tilde{K}_t(x) \Big(\frac{V_t+\Zq_t}{\sigma_w^2}- 1\Big)  \ud t
\end{align*}
where the gain vectors $K(\cdot)$ and $\tilde{K}(\cdot)$ solve the  finite-state space counterpart of the Poisson equations~\eqref{eq:Poisson-eqns}
\begin{align*}
K_t(x) - \bpi_t(x)\sum_{y \in \mathcal X} K_t(y) &= \frac{1}{\sigma_w^2}\bpi_t(x) (h(x) -\barh_t)\\
\tilde{K}_t(x) - \bpi_t(x)\sum_{y \in \mathcal X} \tilde{K}_t(y) &= \frac{1}{\sigma_w^2}\bpi_t(x) (g(x) -\barg_t )
\end{align*}
The general form of the solution is explicitly known:
\begin{align*}
K_t(x) & = \frac{1}{\sigma_w^2}\bpi_t(x)(h(x) - c),\quad 
\tilde{K}_t(x) = \frac{1}{\sigma_w^2}\bpi_t(x)(g(x) - \tilde{c})
\end{align*}
where $c$ and $\tilde{c}$ are constants. The constants are chosen so that $U_t^{x}$ and $\tilde{U}_t^x $ are non-decreasing leading to a well-posed Poisson processes $N^{x,y}(U^x_t)$ and $\tilde{N}^{x,y}(\tilde{U}^x_t)$.  In particular,
\begin{align*}
c = \begin{cases} 
\min_x h(x) ,\quad \text{if}\quad \ud \hat{Z}_t- \barh_t \ud t \geq 0\\
\max_x h(x), \quad \text{else}
\end{cases},\quad
\tilde{c} = \begin{cases} 
\min_x g(x)   ,\quad\text{if}\quad  (\frac{V_t+\Zq_t}{\sigma_w^2}- 1) \geq 0\\
\max_x g(x),  \quad \text{else}
\end{cases}
\end{align*}

\begin{remark}
The FPF for finite state-space Markov chain is proposed in~\cite{yang2016MarkovChain}. It simulates the Wonham filter. Notice that the first line of~\eqref{eq:FPF-finite} has the same structure with the algorithm proposed in~\cite{yang2016MarkovChain} and it is indeed identical when $M=1$. % For multiple agents, additional term is introduced to simulate the effect from the variance of the observations. 
\end{remark}

\section{Simulations}\label{sec:sims}

In this section, we simulate the collective filtering algorithm for
a simple linear-Gaussian system.  There are two objectives: (i) To evaluate
the collective filter described in Prop.~\ref{prop:continuous-LG} as
the number of agents $M$ increases; (ii) To show the convergence of the
estimates using the FPF algorithm as the number of
particles $N\to\infty$.

Comparisons of the collective filtering algorithm are made against the
gold standard of running independent Kalman filters with {\em known}
data association.  It will also be interesting to compare the results
using joint probabilistic data association (JPDA) filter and this is 
planned as part of the continuing research.  

% but those
% results could not be completed in time for the submission.  If the
% paper is accepted, they will be included in the final version of the
% paper.  

The continuous-time system~\eqref{eq:CLG-system} is simulated using
the parameters
\[
A = \begin{pmatrix}
0 & 1\\
-1 & -0.5
\end{pmatrix},\quad H = \begin{pmatrix}
0 & 1
\end{pmatrix}
\]
The process noise covariance $Q = 0.1 I$ and the measurement noise
covariance $\sigma_w^2 = 0.7$.
The initial condition is sampled from a Gaussian prior with
parameters 
\[
m_0 = \begin{pmatrix}
1 \\ 0
\end{pmatrix},\quad \Sigma_0 = \begin{pmatrix}
1 & 0.2 \\ 0.2 & 1
\end{pmatrix}
\]
The sample path for each agent is generated by using an Euler method
of numerical integration with a fixed step size $\Delta t = 0.01$ over
the total simulation time interval $[0,5]$. At  each discrete
time-step, $q$ is approximated as a Gaussian whose mean and variance
are defined as follows:
\begin{align*}
\hat{Z}_t &= \frac{1}{M}\sum_{j=1}^M Z_t^j\\
V_t &= \frac{1}{\Delta t M}\sum_{j=1}^M \big((Z_{t+\Delta t}^j-\hat{Z}_{t+\Delta t})-(Z_{t}^j-\hat{Z}_{t})\big)^2
\end{align*}

The comparison is carried out for the following three filtering algorithms:
\begin{romannum}
\item (KF) This involves simulating $M$ independent Kalman-Bucy
  filters for the $M$ sample paths $\{Z^j_t:0\leq t\leq 5\}$ for
  $j=1,2,\hdots,M$.  The data association is fully known.  
\item (CKF) This involves simulating the mean and the variance of a
  single collective Kalman-Bucy filter using the filtering equations~\eqref{eq:ckf} in Prop.~\ref{prop:continuous-LG}.
\item (FPF) This involves simulating a single FPF~\eqref{eq:FPF}
  with $N$ particles.
\end{romannum}

At the terminal time $T$, KF simulation yields $M$ Gaussians
(posterior distributions for each of the $M$ independent Kalman
filters) whose mean and variance are $m_T^{(j)}$ and $\Sigma_T^{(j)}$, respectively for $j=1,2,\ldots,M$.  We use $m_T^{\text{kf}}$ and $\Sigma_T^{\text{kf}}$ to
denote the mean and the variance of the sum
(mixture) of these $M$ Gaussians. Note the mean and the variance is computed by:
\begin{align*}
m_T^{\text{kf}} &= \frac{1}{M} \sum_{j=1}^M m_T^{(j)}\\
\Sigma_T^{\text{kf}} &= \frac{1}{M} \sum_{j=1}^M \Sigma_T^{(j)} + \big(m_T^{(j)}-m_T^{\text{kf}}\big)\big(m_T^{(j)}-m_T^{\text{kf}}\big)^\top
\end{align*}
The mean and the variance for the CKF is denoted as $m_T^{\text{ckf}}$ and $\Sigma_T^{\text{ckf}}$,
  respectively.  Similarly, $m_T^{\text{fpf}}$ and
  $\Sigma_T^{\text{fpf}}$ are the empirical mean and variance computed
  using the FPF algorithm with $N$ particles.  

% The posterior computed from the collective Kalman-Bucy filter is denoted by ${\cal N}\big(m_T^{\text{ckf}},\Sigma_T^{\text{ckf}}\big)$. Finally, $m_T^{\text{fpf}}$ and $\Sigma_T^{\text{fpf}}$ are the empirical mean and variance from FPF. The distance is mean-squared difference computed from 10,000 simulations, normalized by mean and variance of Kalman filter:

\begin{figure}
\centering
        \centering
        \includegraphics[width=0.45\textwidth]{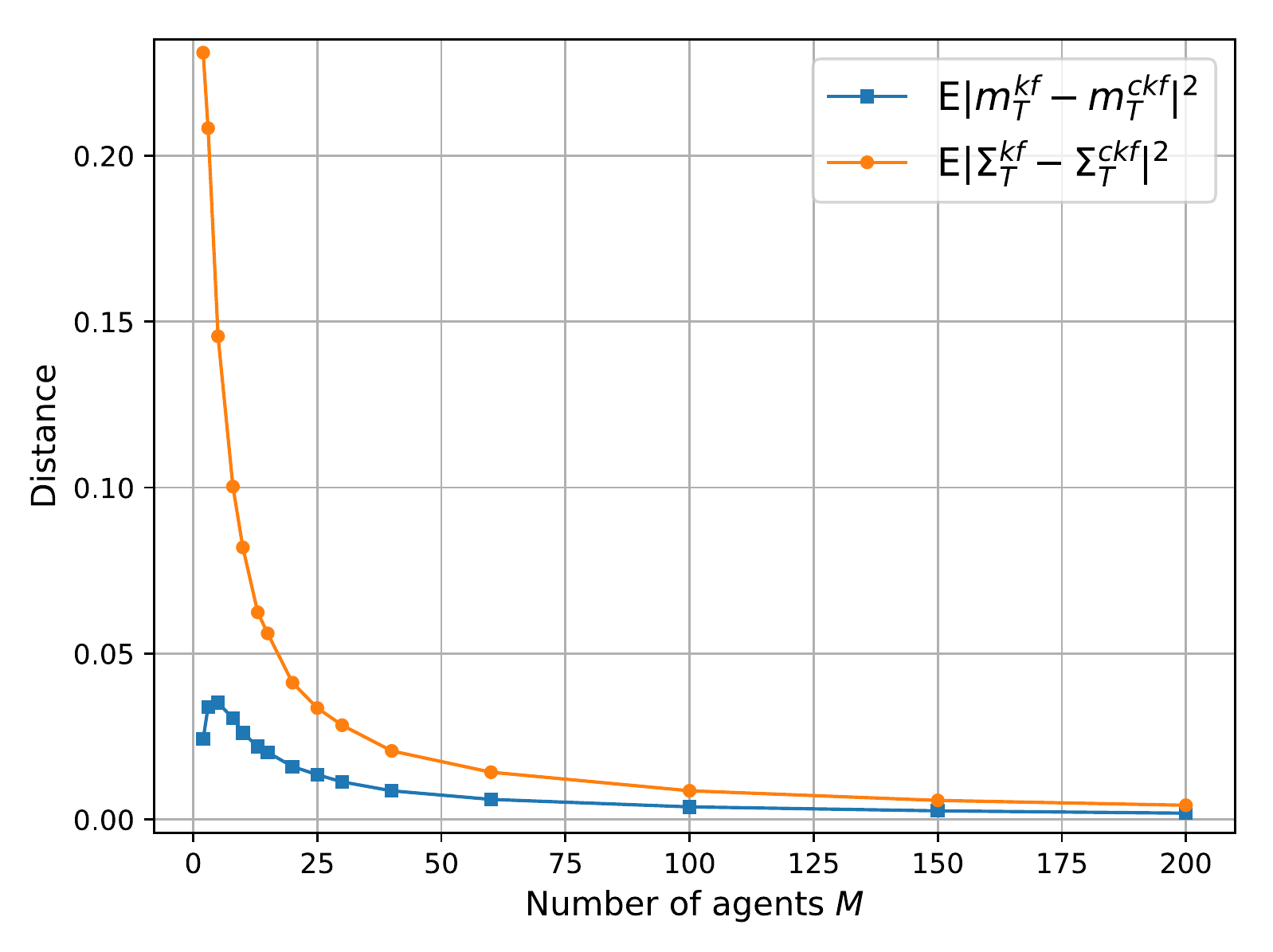}
    %\caption{The distance of mean (blue square) and variance (orange circle) between independent Kalman-Bucy filters with full association and collective Kalman filter described in Proposition~\ref{prop:continuous-LG}.}
        \includegraphics[width=0.45\textwidth]{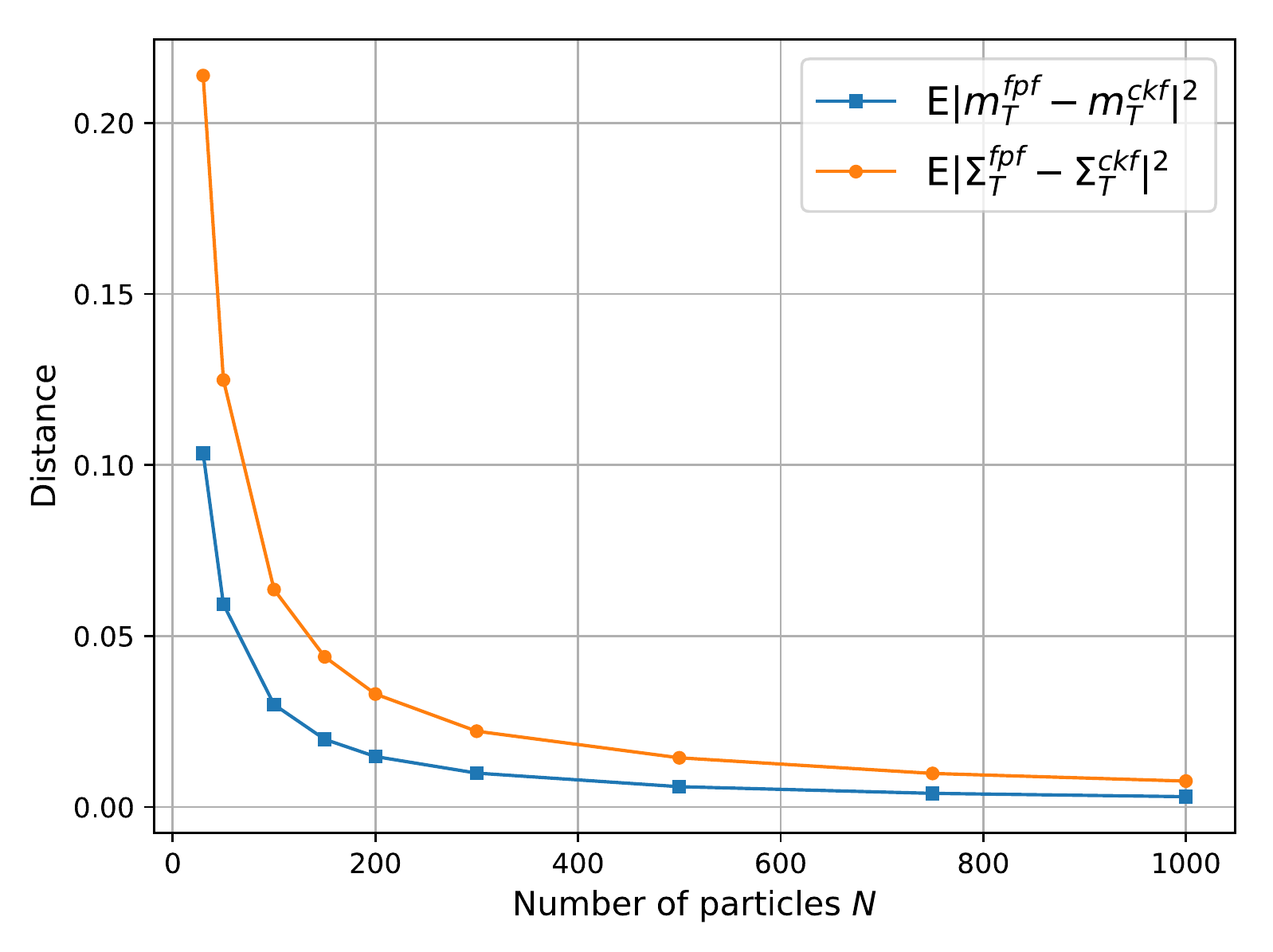}
    %\caption{The distance of mean (blue square) and variance (orange circle) between independent collective Kalman filter (Proposition~\ref{prop:continuous-LG}) and feedback particle filter~\eqref{eq:FPF_impl}.}
\caption{Normalized error for mean (blue circle) and variance (orange circle) with the KF and CKF algorithms.  The KF algorithms were run as $M$ independent Kalman filters with fully known data associations.}
\label{fig:change-M}
\caption{Normalized error for mean (blue circle) and variance (orange circle) with the CKF and FPF algorithms.  The number of agents is fixed to $M=30$ for this simulation.}
\label{fig:change-N}
\end{figure}

Figure~\ref{fig:change-M} depicts the normalized difference for the
mean and the variance between (KF) and (CKF), as the number of agents
increase from $M=2$ to $M=200$.  It is observed that both the
differences converge to zero as the number of agents increases.  Though omitted in the plot, for the $M=1$ case, these do match exactly.

Figure~\ref{fig:change-N} depicts the distance between normalized difference for the
mean and the variance between (CKF) and (FPF).  In this simulation,
$M=30$ is fixed and $N$ varies from 30 to 1000.  The plots show that
the difference converges to zero as $N$ increases.  Therefore, the FPF
is able to provide a solution to the collective inference problem.

% \section{Conclusion and Discussions}

% In this paper, recursive update to obtain optimal estimate for the collective filtering problem is presented. The optimization objective~\eqref{eq:KL-optimization} is a special case of the problem introduced in~\cite{singh2020inference} with the size of time window is 1. 
% The algorithm is naturally extended to continuous time system, where the standard nonlinear filter is a special case of single agent ($M=1$). For continuous-time system, the FPF algorithm is proposed to numerically approximate the collective filter.

% (Additional discussion on utility of this work?)

\appendix

\section{Proofs of the propositions}

\subsection{Proof of proposition~\ref{prop:one-step}}\label{apdx:one-step-estimator}
Eq.~\eqref{eq:one-step-estimator-a} for $\pi_{t+1\mid t}$ is the usual prediction step, and the nominal distribution $\sP(x,z) = \pi_{t+1\mid t}(x) o(z\mid x)$. Therefore the optimization problem becomes
\[
\min_{\tsP} \sum_{x,z}\tsP(x,z)\log\frac{\tsP(x,z)}{\pi_{t+1\mid t}(x) o(z\mid x)}
\]
and the constraint is that the marginal of $\tsP$ on $Z_{t+1}$ must be $q_{t+1}$.

\medskip

Hence a Lagrange multiplier $\lambda(z)$ is introduced and the objective function becomes
\[
\sum_{x,z}\tsP(x,z)\log\frac{\tsP(x,z)}{\pi_{t+1\mid t}(x) o(z\mid x)}+ \sum_z \lambda(z)\big(\sum_x\tsP(x,z)-q_{t+1}(z)\big)
\] 
Differentiate with respect to $\tsP(x,z)$ yields
$$
\log\frac{\tsP(x,z)}{\pi_{t+1\mid t}(x) o(z\mid x)} + 1 + \lambda(z) = 0
$$
The solution is
$$
\tsP(x,z) = \pi_{t+1\mid t}(x) o(z\mid x)\exp\big(-1-\lambda(z)\big)
$$
It is substituted to the constraint;
$$
\sum_x \pi_{t+1\mid t}(x) o(z\mid x)\exp\big(-1-\lambda(z)\big) - q_{t+1}(z) = 0
$$
Therefore,
$$
\exp\big(-1-\lambda(z)\big) = \frac{q_{t+1}(z)}{\sum_x \pi_{t+1\mid t}(x) o(z\mid x)}
$$
Denote the denominator by $\xi_{t+1}(z) = \sum_x \pi_{t+1\mid t}(x) o(z\mid x)$ and collect the result to conclude:
$$
\pi_{t+1\mid t+1}(x) = \sum_z \tsP(x,z) = \sum_z \frac{\pi_{t+1\mid t}(x)o(z\mid x)q_{t+1}(z)}{\xi_{t+1}(z)}
$$

\subsection{Proof of proposition~\ref{prop:continuous-time}}\label{apdx:continuous-time}

 The continuous-time limit of the first update step~\eqref{eq:one-step-estimator-a} is straightforward. 
 \begin{equation*}\label{eq:forward-step}
 \frac{\ud}{\ud t}\pi_t(x) = \clA^\dagger \pi_t(x)
 \end{equation*}
In order to derive the continuous-time limit of the update step~\eqref{eq:one-step-estimator-b}, 
%  recall the observation model
% \begin{equation*}
% 	\ud Z_t  = h(X_t)\ud t + \sigma_w \ud W_t
% \end{equation*}
% concluding the observation likelihood
the likelihood function is expressed as a Gaussian
\[
o(Z_{t+\Delta t} - Z_t|X_t) \propto \exp\Big(-\frac{|Z_{t+\Delta t} - Z_t - \int_t^{t+\Delta t} h(X_s)\ud s |^2}{2\sigma_w^2\Delta t}\Big)
\]
The update formula for $t \to t+\Delta t$ is 
\begin{align*}
\pi_{t+\Delta t}(x) = \clA^\dagger \pi_t(x) \Delta t + \int_{\Delta Z} \frac{\exp(-\frac{|\Delta Z - h(x)\Delta t + O(\Delta t^2)|^2}{2\sigma_w^2\Delta t}) \pi_t(x)}{\xi_t(\Delta Z)} \ud q_t(\Delta Z) + O(\Delta t^2)
\end{align*}
where the approximation $\int_t^{t+\Delta t} h(X_s)\ud s = h(X_t)\Delta t + O(\Delta t^2)$ is used. The denominator
\begin{equation*}
\xi_t(\Delta Z) = \int \exp\big(-\frac{|\Delta Z - h(x)\Delta t + O(\Delta t^2)|^2}{2\sigma_w^2\Delta t}\big) \pi_t(x)\ud x
\end{equation*}
Note that the term  $\exp(-\frac{\Delta Z^2}{2\sigma_w^2\Delta t})$ appears both in numerator and denominator and cancel.
For the rest of the terms, we use the expansion keeping the terms up to second-order for $\Delta Z$ and first-order for $\Delta t$. Since $\exp(x) = 1 + x + \frac{1}{2}x^2 + O(x^3)$, the numerator inside the integral becomes:
\begin{align*}
\exp\Big(-\frac{ h(x)^2\Delta t -2h(x)\Delta Z + O(\Delta t^2)}{2\sigma_w^2}\Big) &\pi_t(x) \\
= (1-\frac{1}{2\sigma_w^2}h(x)^2 \Delta t + \frac{1}{\sigma_w^2}h(x) &\Delta Z + \frac{1}{2\sigma_w^4}h(x)^2\Delta Z^2 + O(\Delta t^2,\Delta Z^3))\pi_t(x)
\end{align*}
Similarly the denominator is also expressed by:
\begin{equation*}
\int \big(1-\frac{1}{2\sigma_w^2}h(x)^2 \Delta t + \frac{1}{\sigma_w^2}h(x) \Delta Z + \frac{1}{2\sigma_w^4}h(x)^2\Delta Z^2 + O(\Delta t^2,\Delta Z^3)\big)\pi_t(x) \ud x
\end{equation*}
We use $\frac{1+x}{1+y} = 1+x-y +y^2 -xy + O(xy^2,y^3)$ to express the ratio, and then it simplifies to
\begin{align*}
\pi_{t+\Delta t}(x) &= \clA^\dagger \pi_t(x) \Delta t + \int_{\Delta Z} \frac{1}{\sigma_w^2}\pi_t(x)(h(x) - \hat{h}_t)(\Delta Z - \hat{h}_t\Delta t) \ud q_t(\Delta Z) \\&+ \int_{\Delta Z} \frac{1}{\sigma_w^2}\pi_t(x)(g(x) - \hat{g}_t)(\frac{1}{\sigma_w^2}\Delta Z^2 - \Delta t)\ud  q_t(\Delta Z) \\
&+  O(\Delta t^2,\Delta Z^3)
\end{align*}
%\begin{align*}
%\pi_{t+\Delta t}(x) &= \int_{\Delta Z} \frac{\exp(-\frac{ h(x)^2\Delta t -2h(x)\Delta Z + O(\Delta t^2)}{2\sigma_w^2}) \pi_t(x)}{ \sum_{x'}\exp(-\frac{ h(x')^2\Delta t -2h(x')\Delta Z + O(\Delta t^2)}{2\sigma_w^2})  \pi_t(x')} q_t(\Delta Z)\\
%&=\int_{\Delta Z} \frac{(1-\frac{1}{2\sigma_w^2}h(x)^2 \Delta t + \frac{1}{\sigma_w^2}h(x) \Delta Z + \frac{1}{2\sigma_w^4}h(x)^2\Delta Z^2 + O(\Delta t^2,\Delta Z^3))\pi_t(x)}{ \sum_{x'}(1-\frac{1}{2\sigma_w^2}h(x')^2 \Delta t +\frac{1}{\sigma_w^2} h(x') \Delta Z + \frac{1}{2\sigma_w^4}h(x')^2\Delta Z^2 + O(\Delta t^2,\Delta Z^3))\pi_t(x')} q_t(\Delta Z)\\
%&= \int_{\Delta Z} \pi_t(x)\bigg[ 1-\frac{1}{2\sigma_w^2}h(x)^2 \Delta t + \frac{1}{\sigma_w^2}h(x) \Delta Z + \frac{1}{2\sigma_w^4}h(x)^2\Delta Z^2  \\
% &\quad \quad \quad -  \sum_{x'} (-\frac{1}{2\sigma_w^2}h(x')^2 \Delta t + \frac{1}{\sigma_w^2}h(x') \Delta Z + \frac{1}{2\sigma_w^4}h(x')^2\Delta Z^2)\pi_t(x') +  (\sum_{x'}  \frac{1}{\sigma_w^2}h(x')\pi_t(x') \Delta Z)^2 \\
% &  \quad \quad \quad -  \frac{1}{\sigma_w^4}h(x) \sum_{x'} h(x')\pi_t(x') \Delta Z^2  \bigg]q_t(\Delta Z) +  O(\Delta t^2,\Delta Z^3)\\
% &=\pi_t(x)  + \int_{\Delta Z} \frac{1}{\sigma_w^2}\pi_t(x)(h(x) - \hat{h}_t)(\Delta Z - \hat{h}_t\Delta t) q_t(\Delta Z) \\&+ \int_{\Delta Z} \frac{1}{\sigma_w^2}\pi_t(x)(g(x) - \hat{g}_t  - h(x)\hat{h}_t + \hat{h}_t^2)(\frac{1}{\sigma_w^2}\Delta Z^2 - \Delta t)q_t(\Delta Z) +  O(\Delta t^2,\Delta Z^3)
%\end{align*}
%where we used $\exp(x) = 1 + x + \frac{1}{2}x^2 + O(x^3)$ and $\frac{1+x}{1+y} = 1+x-y +y^2 -xy + O(xy^2,y^3)$. 
Dividing by $\Delta t$ and taking the limit as $\Delta t \to 0$ concludes the result.

\subsection{Proof of proposition~\ref{prop:continuous-LG}}\label{apdx:pf-continuous-LG}

For linear-Gaussian example, the mean and the variance fully
characterize the distribution. Thus, we repeat the procedure
in~\ref{sec:ctmc} for the linear-Gaussian case.  Although the linear
Gaussian is a special case of the general Euclidean result, we provide
here a proof as a continuous limit of a discrete-time model (which is
of independent interest).  The results are stated and proved in somewhat more
general settings with vector-valued observations.

\subsubsection{Discrete-time linear-Gaussian problem} The model is
\begin{subequations}\label{eq:DLG-system}
\begin{align}
X_{t+1} &= AX_t + B_t,\quad X_0\sim {\cal N}(m_0,\Sigma_0)\\
Z_t &= HX_t + W_t
\end{align}
\end{subequations}
where $\{B_t\}_{t\geq 0}$, $\{W_t\}_{t\geq 0}$ are mutually
independent i.i.d. Gaussian
random variables with zero mean and variance $Q$ and $R$,
respectively, and also assumed to be independent of $X_0$. 
The observation $q_t$ is assumed to be a Gaussian with mean $\hat{Z}_t$ and variance $V_t$. The discrete-time update~\eqref{eq:one-step-estimator} for the mean and the variance of linear-Gaussian problem is illustrated in the following proposition.

\begin{proposition}\label{prop:discrete-LG}
Consider the collective filtering problem for the discrete-time linear
Gaussian model~\eqref{eq:DLG-system}.  Suppose  $q_t = {\cal
  N}( \hat{Z}_t, V_t)$ and $\pi_{t\mid t}= {\cal
  N}( m_{t\mid t}, \Sigma_{t\mid t})$ are both Gaussian.  Then $\pi_{t+1\mid t}$
and $\pi_{t+1\mid t+1}$ are also Gaussian whose mean and variance evolve
according to the following recursion:
\begin{subequations}
\begin{align}
m_{t+1\mid t} &= Am_{t\mid t}\\
\Sigma_{t+1\mid t} &= A\Sigma_{t\mid t}A^\top + Q\\
K_{t+1} &= \Sigma_{t+1\mid t} H^\top(H\Sigma_{t+1\mid t} H^\top + R)^{-1}\\
m_{t+1\mid t+1} &=m_{t+1\mid t} +K_{t+1} (\hat{Z}_{t+1}-Hm_{t+1\mid t}) \label{eq:dlg-mean}\\
\Sigma_{t+1\mid t+1} &=\Sigma_{t+1\mid t} - K_{t+1}(H\Sigma_{t+1\mid t} H^\top + R - V_{t+1})K_{t+1}^{-1} \label{eq:dlg-var}
\end{align}
\end{subequations}
\end{proposition}

{\it (proof):} Recall the one-step estimator~\eqref{eq:one-step-estimator} in
continuous state-space settings:
\begin{align*}
\pi_{t+1\mid t+1}(x) &\propto \int \frac{o(z\mid x)\pi_{t+1\mid t}(x)}{\xi(z)}q_{t+1}(z)\ud z\\
&= \pi_{t+1\mid t}(x) \int \frac{o(z\mid x)}{\xi(z)}q_{t+1}(z)\ud z\\
\xi(z) &=\int o(z\mid x)\pi_{t+1\mid t}(x)\ud x
\end{align*}
where the probability density is involved instead of probability mass function.
Note $\xi(z)$ is the pdf of a Gaussian with mean $Hm_{t+1\mid t}$ and variance $H\Sigma_{t+1\mid t}H^\top +R$, and therefore
\begin{align*}
o(z\mid x) &\propto \exp\Big(-\frac{1}{2} (z-Hx)^\top R^{-1}(z-Hx)\Big)\\
\xi(z) \propto &
\exp\Big(-\frac{1}{2} (z-Hm_{t+1\mid t})^\top(H\Sigma_{t+1\mid t}H^\top + R)^{-1}(z-Hm_{t+1\mid t})\Big)\\
q_{t+1}(z) &\propto \exp\Big(-\frac{1}{2} (z-\hat{Z}_{t+1})^\top V_{t+1}^{-1} (z-\hat{Z}_{t+1})\Big)
\end{align*}
Therefore, the integrand $\frac{o(z\mid x)q_{t+1}(z)}{\xi(z)} \propto \exp\big(-\frac{1}{2} E_1\big)$ where
\begin{align*}
E_1 &= (z-Hx)^\top R^{-1}(z-Hx)+ (z-\hat{Z}_{t+1})^\top V_{t+1}^{-1} (z-\hat{Z}_{t+1})\\
&\quad- (z-Hm_{t+1\mid t})^\top(H\Sigma_{t+1\mid t}H^\top + R)^{-1}(z-Hm_{t+1\mid t}) \\
&= (z-c_0)^\top C_0(z-c_0) + (Hx-\hat{Z}_{t+1})^\top C_1 (Hx-\hat{Z}_{t+1})\\
&\quad + (Hm_{t+1\mid t}-\hat{Z}_{t+1})^\top C_2(Hm_{t+1\mid t}-\hat{Z}_{t+1})\\
&\quad + (Hx-Hm_{t+1\mid t})^\top C_3 (Hx-Hm_{t+1\mid t}) 
\end{align*}
where
\begin{align*}
C_1 &= R^{-1}\big(R^{-1}+V_{t+1}^{-1}-(H\Sigma_{t+1\mid t}H^\top + R)^{-1}\big)^{-1}V_{t+1}^{-1}\\
C_2 &= -V_{t+1}^{-1}\big(R^{-1}+V_{t+1}^{-1}-(H\Sigma_{t+1\mid t}H^\top + R)^{-1}\big)^{-1}(H\Sigma_{t+1\mid t}H^\top + R)^{-1}\\
C_3 &= -R^{-1}\big(R^{-1}+V_{t+1}^{-1}-(H\Sigma_{t+1\mid t}H^\top + R)^{-1}\big)^{-1}(H\Sigma_{t+1\mid t}H^\top + R)^{-1}
\end{align*}
%\begin{align*}
%C_1 &= \frac{H\Sigma_{t+1\mid t}H^\top + R}{R(H\Sigma_{t+1\mid t}H^\top + R)+H\Sigma_{t+1\mid t}H^\top V_{t+1}}\\
%C_2 &= \frac{-R}{R(H\Sigma_{t+1\mid t}H^\top + R)+H\Sigma_{t+1\mid t}H^\top V_{t+1}}\\
%C_3 &= \frac{-V_{t+1}}{R(H\Sigma_{t+1\mid t}H^\top + R)+H\Sigma_{t+1\mid t}H^\top V_{t+1}}
%\end{align*}
%Here, $\int \exp\big(-\frac{1}{2} (z-c_0)^\top C_0 (z-c_0)\big)\ud z$
%becomes a constant, and all other terms are pulled out from the
%integral. 
%JINKIM -- Do not understand the above
Also, $\pi_{t+1\mid t+1}(x) \propto \exp\big(-\frac{1}{2} E_2\big)$ where
\begin{align*}
E_2&= (x-m_{t+1\mid t})^\top \Sigma_{t+1\mid t}^{-1}(x-m_{t+1\mid t}) + (Hx-\hat{Z}_{t+1})^\top C_1 (Hx-\hat{Z}_{t+1})\\
&\quad + (Hx-Hm_{t+1\mid t})^\top C_3 (Hx-Hm_{t+1\mid t})\\
&=x^\top\big(\Sigma_{t+1\mid t}^{-1}+H^\top(C_1+C_3) H\big)x\\
&\quad -2x^\top\big(\Sigma_{t+1\mid t}^{-1}m_{t+1\mid t} + H^\top C_1 \hat{Z}_{t+1} + H^\top C_3 H m_{t+1\mid t}\big) + \text{(const.)}
\end{align*}
Therefore, $\pi_{t+1\mid t+1}$ is a Gaussian pdf with mean and variance is given by:
\begin{align}
m_{t+1\mid t+1} &=\Sigma_{t+1\mid t+1}\big(\Sigma_{t+1\mid t}^{-1}m_{t+1\mid t} + H^\top C_1 \hat{Z}_{t+1} + H^\top C_3 H m_{t+1\mid t}\big)\label{eq:mean-interm}\\
\Sigma_{t+1\mid t+1} &= \big(\Sigma_{t+1\mid t}^{-1}+H^\top(C_1+C_3) H\big)^{-1}\label{eq:var-interm}
\end{align}
By the matrix inversion lemma, the variance formula becomes
\begin{align*}
\Sigma_{t+1\mid t+1} &= \Sigma_{t+1\mid t} - \Sigma_{t+1\mid t} H^\top \big((C_1+C_3)^{-1}+H\Sigma_{t+1\mid t} H^\top\big)^{-1} H\Sigma_{t+1\mid t}
\end{align*}
%We use a fact that for any matrix $F$ and $G$,
%$$
%F^{-1} + G^{-1}= F^{-1}GG^{-1} + F^{-1}FG^{-1}=F^{-1}(G+F)G^{-1}
%$$
%provided that the inverses exist.
Observe that
\begin{align*}
C_1+C_3 &= R^{-1}\big(R^{-1}+V_{t+1}^{-1}-(H\Sigma_{t+1\mid t}H^\top + R)^{-1}\big)^{-1}\\
&\quad\cdot \big(V_{t+1}^{-1}-(H\Sigma_{t+1\mid t}H^\top + R)^{-1}\big)\\
&= \big(R+(V_{t+1}^{-1}-(H\Sigma_{t+1\mid t}H^\top + R)^{-1})^{-1}\big)^{-1}
\end{align*}
Therefore
\begin{align*}
(C_1+&C_3)^{-1}+H\Sigma_{t+1\mid t} H^\top\\
&= R - ((H\Sigma_{t+1\mid t} H^\top + R)^{-1}-V_{t+1}^{-1})^{-1} + H\Sigma_{t+1\mid t} H^\top\\
&=- V_{t+1}(H\Sigma_{t+1\mid t} H^\top + R-V_{t+1})^{-1}(H\Sigma_{t+1\mid t} H^\top + R)\\
&\quad+(H\Sigma_{t+1\mid t} H^\top + R)\\
&=(H\Sigma_{t+1\mid t} H^\top + R)(H\Sigma_{t+1\mid t} H^\top + R-V_{t+1})^{-1}\\
&\quad \cdot(H\Sigma_{t+1\mid t} H^\top + R)
\end{align*}
Substituting back to the variance equation~\eqref{eq:var-interm}, the Riccati equation~\eqref{eq:dlg-var} is obtained. It is substituted to~\eqref{eq:mean-interm} to obtain~\eqref{eq:dlg-mean}.

\subsubsection{Proof of proposition~\ref{prop:continuous-LG}} The previous proposition is extended to continuous-time problem by
considering suitable limits.

\def\dtstep{{\Delta t}}

Consider the continuous-time system~\eqref{eq:CLG-system} with a
discrete time-step $\dtstep$,
\begin{align*}
X_{t+\dtstep} &= (I+A\dtstep) X_t + \Delta{B}_t\\
Z_{t+\dtstep} &= Z_t+H X_t\dtstep + \Delta{W}_t
\end{align*}
where $\Delta{B}_t$ and $\Delta{W}_t$ are normal random variables with variance $Q\dtstep$, $\sigma_w^2 I \dtstep$, respectively. $Z_{t+\dtstep} - Z_t$ is assumed to be a normal random variable and its mean and variance are $\Delta \hat{Z}_t:=\hat{Z}_{t+\dtstep}=\hat{Z}_t$, and $V_t\dtstep$ respectively.

By the Proposition~\ref{prop:discrete-LG}, the prediction step is:
\begin{align*}
m_{t+\dtstep\mid t} &= (I+A\dtstep)m_{t}\\
\Sigma_{t+\dtstep\mid t} &= (I+A\dtstep)\Sigma_{t}(I+A\dtstep)^\top + Q\dtstep
\end{align*}
For the estimation step, omit the higher order terms such as $H \Sigma_{t+1\mid t} H^\top \dtstep^2$ to simplify the equation, and then we have
\begin{align*}
m_{t+\dtstep\mid t+\dtstep} &=m_{t+\dtstep\mid t} -\Sigma_{t+\dtstep\mid t} H^\top R^{-1} (\Delta \hat{Z}_{t} - Hm_{t+\dtstep\mid t}\dtstep) \\
\Sigma_{t+\dtstep\mid t+\dtstep} &= \Sigma_{t+\dtstep\mid t}- \Sigma_{t+\dtstep\mid t} H^\top R^{-1}(R - V_t)R^{-1} H \Sigma_{t+\dtstep\mid t}\dtstep
\end{align*}
%\begin{align*}
%m_{t+\dtstep\mid t+\dtstep} &=m_{t+\dtstep\mid t} +\Sigma_{t+\dtstep\mid t} H^\top\dtstep(H\Sigma_{t+\dtstep\mid t} H^\top\dtstep^2 + R\dtstep)^{-1} (\Delta \hat{Z}_{t} - Hm_{t+\dtstep\mid t}\dtstep)\\
%\Sigma_{t+\dtstep\mid t+\dtstep} &= \Sigma_{t+\dtstep\mid t} - \Sigma_{t+\dtstep\mid t} H^\top \dtstep\Big(\frac{(H\Sigma_{t+\dtstep\mid t} H^\top\dtstep^2 + R\dtstep)^2}{H\Sigma_{t+\dtstep\mid t} H^\top\dtstep^2 + R\dtstep - V_t \dtstep}\Big)^{-1} H\dtstep\Sigma_{t+\dtstep\mid t}
%\end{align*}
up to $o(\dtstep)$ error. 
Substitute $m_{t+\dtstep\mid t}$ and $\Sigma_{t+\dtstep\mid t}$ to the equation, and ignoring higher order terms,
\begin{align*}
&m_{t+\dtstep} - m_{t} = Am_{t}\dtstep+ \Sigma_{t} H^\top R^{-1}(\Delta \hat{Z}_t - Hm_{t}\dtstep) \\
&\Sigma_{t+\dtstep} - \Sigma_{t} =\Big(A\Sigma_{t} + \Sigma_{t}A^\top + Q - \Sigma_{t} H^\top \big(R^{-1}-R^{-1}V_tR^{-1}\big)H\Sigma_{t} \Big)\dtstep 
\end{align*}
The differential formula is obtained by letting $\dtstep \to 0$, and $R = \sigma_w^2$.

\subsection{Proof of proposition~\ref{prop:FPF-Euc}}\label{apdx:pf-FPF-Euc}

In order to check if the FPF update law gives the required distribution, we express the Fokker-Planck equation for $X^i_t$. The It\^o form of the Stratonovich sde~\eqref{eq:FPF-euc} is
\begin{align*}
\ud X_t^i  = & a(X_t^i)\ud t + \sigma(X_t^i) \ud B_t^i + K_t(X_t^i)\big(\ud \hat{Z}_t - \barh_t \ud t \big)  - \alpha_t K_t(X_t^i) \big(h(X_t^i) -\barh_t\big)\ud t\nonumber\\
&+\frac{1}{2} \sum_{n=1}^d K^n_t (X^i_t)\frac{\partial K_t}{\partial x_n} (X^i_t) \Zq_t \ud t+ \Big(\frac{V_t+\Zq_t}{\sigma_w^2}- 1\big)u_t(X_t^i) \ud t
\end{align*}
The additional term $\frac{1}{2} \sum_{n=1}^d K^n_t (X^i_t)\frac{\partial K_t}{\partial x_n} (X^i_t) \Zq_t \ud t$ appears due to change from Stratonovich to It\^o form, where $K_t^n(x)$ denotes the $n$-th component of the gain vector $K_t(x)=(K^1_t(x),\ldots, K^d_t(x)) \in \Re^d$, and $\hat{V}_t \ud t$ is the quadratic variation of the stochastic process $\hat{Z}_t$. The corresponding Fokker-Planck equation is  
\begin{align*}
\ud \pi_t =& \mathcal A^\dagger  \pi_t \ud t- \nabla \cdot(\pi_t  K_t)  ( \ud \hat{Z}_t - \barh_t \ud t)  +   \frac{1}{2}  \sum_{n,m=1}^d \frac{\partial^2}{\partial x_n \partial x_m}( K^n_t  K^m_t \pi_t) \Zq_t \ud t   \\
&+ \alpha_t \nabla \cdot ( \pi_t  K_t (h-\barh_t)) \ud t- \frac{1}{2} \sum_{n,m=1}^d \frac{\partial}{\partial x_n} (\pi_t  K^m_t \frac{\partial}{\partial x_m}  K^n_t) \Zq_t \ud t \\
&- \Big(\frac{V_t+\Zq_t}{\sigma_w^2}- 1\Big)\nabla \cdot(\pi_t  u_t)  \ud t 
%  =& \frac{1}{\sigma_w^2}\pi_t(h-\barh_t)( \ud \hat{Z}_t - \barh_t \ud t) + \frac{1}{2}  \sum_{i,j} \frac{\partial}{\partial x_i}( K_i \frac{\partial}{\partial x_j} ( K_j \pi_t)) \Zq_t \ud t  \\
%&+  \alpha_t \nabla \cdot ( \pi_t  K_t (h-\barh_t)) \ud t +   \Big(\frac{V_t+\Zq_t}{\sigma_w^2}- 1\Big) \big(\frac{1}{2}\nabla \cdot(\pi_t  K_t (h-\barh_t)) + \frac{1}{\sigma_w^2}\pi_t(g-\barg_t)\big)\\
%= &\frac{1}{\sigma_w^2}\pi_t(h-\barh_t)( \ud \hat{Z}_t - \barh_t \ud t)  + \frac{1}{2} \nabla \cdot(  K_t \nabla \cdot(\pi_t  K_t)) \Zq_t \ud t +  \alpha_t \nabla \cdot ( \pi_t  K_t (h-\barh_t)) \ud t\\&+  \Big(\frac{V_t+\Zq_t}{\sigma_w^2}- 1\Big) \big(\frac{1}{2}\nabla \cdot(\pi_t  K_t (h-\barh_t)) + \frac{1}{\sigma_w^2}\pi_t(g-\barg_t)\big)\ud t\\
%= &\frac{1}{\sigma_w^2}\pi_t(h-\barh_t)( \ud \hat{Z}_t - \barh_t \ud t)  - \frac{1}{2\sigma_w^2} \nabla \cdot( \pi_t  K_t (h-\barh_t)) \Zq_t \ud t +  \alpha_t \nabla \cdot ( \pi_t  K_t (h-\barh_t)) \ud t\\&+  \frac{1}{2} \Big(\frac{V_t+\Zq_t}{\sigma_w^2}- 1\Big) \nabla \cdot(\pi_t  K_t (h-\barh_t))  \ud t+ \frac{1}{\sigma_w^2}\Big(\frac{V_t+\Zq_t}{\sigma_w^2}- 1\Big)\pi_t(g-\barg_t) \ud t\\
%%    = &\frac{1}{\sigma_w^2}\pi_t(h-\barh_t)( \ud \hat{Z}_t - \barh_t \ud t)  +  \Big(\frac{V_t+\Zq_t}{\sigma_w^2}- 1\Big) \nabla \cdot(  K_t \pi_t \barh_t)  \ud t+ \frac{1}{\sigma_w^2}\Big(\frac{V_t+\Zq_t}{\sigma_w^2}- 1\Big)\pi_t(g-\barg_t)  \\
%= & \frac{1}{\sigma_w^2}\pi_t(h-\barh_t)( \ud \hat{Z}_t - \barh_t \ud t)  + \frac{1}{\sigma_w^2} \Big(\frac{V_t+\Zq_t}{\sigma_w^2}- 1\Big) \pi_t(g-\barg_t ) \ud t
\end{align*}
Upon using the Poisson equations~\eqref{eq:Poisson-eqns} and collecting the terms
\begin{align*}
  \ud \pi_t =
%  & \ud t- \nabla \cdot(\pi_t  K_t)  ( \ud \hat{Z}_t - \barh_t \ud t)  +   \frac{1}{2}  \sum_{n,m} \frac{\partial^2}{\partial x_n \partial x_m}( K^n_t  K^m_t \pi_t) \Zq_t \ud t   \\
% &- \frac{1}{2} \sum_{i,j} \frac{\partial}{\partial x_i} (\pi_t  K_j \frac{\partial}{\partial x_j}  K_i) \Zq_t \ud t+ \alpha_t \nabla \cdot ( \pi_t  K_t (h-\barh_t)) \ud t \\
%&- \Big(\frac{V_t+\Zq_t}{\sigma_w^2}- 1\Big)\nabla \cdot(\pi_t  u_t)  \ud t \\
%  =& \frac{1}{\sigma_w^2}\pi_t(h-\barh_t)( \ud \hat{Z}_t - \barh_t \ud t) + \frac{1}{2}  \sum_{i,j} \frac{\partial}{\partial x_i}( K_i \frac{\partial}{\partial x_j} ( K_j \pi_t)) \Zq_t \ud t  \\
%&+  \alpha_t \nabla \cdot ( \pi_t  K_t (h-\barh_t)) \ud t +   \Big(\frac{V_t+\Zq_t}{\sigma_w^2}- 1\Big) \big(\frac{1}{2}\nabla \cdot(\pi_t  K_t (h-\barh_t)) + \frac{1}{\sigma_w^2}\pi_t(g-\barg_t)\big)\\
  & \mathcal A^\dagger  \pi_t \ud t +\frac{1}{\sigma_w^2}\pi_t(h-\barh_t)( \ud \hat{Z}_t - \barh_t \ud t)  + \frac{1}{2} \nabla \cdot(  K_t \nabla \cdot(\pi_t  K_t)) \Zq_t \ud t \\
&+  \alpha_t \nabla \cdot ( \pi_t  K_t (h-\barh_t)) \ud t +  \Big(\frac{V_t+\Zq_t}{\sigma_w^2}- 1\Big) \big(\frac{1}{2}\nabla \cdot(\pi_t  K_t (h-\barh_t)) \\
&+ \frac{1}{\sigma_w^2}\pi_t(g-\barg_t)\big)\ud t
  \\
   = & \mathcal A^\dagger  \pi_t \ud t+\frac{1}{\sigma_w^2}\pi_t(h-\barh_t)( \ud \hat{Z}_t - \barh_t \ud t)  - \frac{1}{2\sigma_w^2} \nabla \cdot( \pi_t  K_t (h-\barh_t)) \Zq_t \ud t \\&+  \alpha_t \nabla \cdot ( \pi_t  K_t (h-\barh_t)) \ud t+  \frac{1}{2} \Big(\frac{V_t+\Zq_t}{\sigma_w^2}- 1\Big) \nabla \cdot(\pi_t  K_t (h-\barh_t))  \ud t\\&+ \frac{1}{\sigma_w^2}\Big(\frac{V_t+\Zq_t}{\sigma_w^2}- 1\Big)\pi_t(g-\barg_t) \ud t\\
%    = &\frac{1}{\sigma_w^2}\pi_t(h-\barh_t)( \ud \hat{Z}_t - \barh_t \ud t)  +  \Big(\frac{V_t+\Zq_t}{\sigma_w^2}- 1\Big) \nabla \cdot(  K_t \pi_t \barh_t)  \ud t+ \frac{1}{\sigma_w^2}\Big(\frac{V_t+\Zq_t}{\sigma_w^2}- 1\Big)\pi_t(g-\barg_t)  \\
     = & \mathcal A^\dagger  \pi_t \ud t+ \frac{1}{\sigma_w^2}\pi_t(h-\barh_t)( \ud \hat{Z}_t - \barh_t \ud t)  + \frac{1}{\sigma_w^2} \Big(\frac{V_t+\Zq_t}{\sigma_w^2}- 1\Big) \pi_t(g-\barg_t ) \ud t
  \end{align*}
  concluding the update law~\eqref{eq:cont-time-update} for the collective filter.  

\bibliographystyle{AIMS}
\bibliography{duality,filter-stability-observability,jpda}

%%%%%%%%%%%%%%%%%%%%%%%%%%%%%%%%%%%%%%%%%%%%%%%%%%%%%%%%%%%%%%%%%%%%%%%%%%%%%%

\medskip
% The data information below will be filled by AIMS editorial staff
Received xxxx 20xx; revised xxxx 20xx.
\medskip

\end{document}

%% file: _definitions.tex
\def\qed{\hspace*{\fill}~\IEEEQED\par\endtrivlist\unskip}%\mbox{\rule[0pt]{1.3ex}{1.3ex}}}

\def\Re{\mathbb{R}}

\def\notes#1{\marginpar{\tiny #1}\typeout{Notes!
Notes!
Notes!
}}
\renewcommand{\notes}[1]{\typeout{notes!}}

\def\Re{\field{R}}

\def\clZ{{\cal Z}}

\def\E{{\sf E}}

\def\clX{{\cal X}}

\def\IEEEQEDclosed{\mbox{\rule[0pt]{1.3ex}{1.3ex}}}
\def\qed{\nobreak\hfill\IEEEQEDclosed}

\def\clZ{{\cal Z}}

\def\beq{\begin{eqnarray}} 
\def\bc{\begin{center}} 
\def\be{\begin{enumerate}}
\def\bi{\begin{itemize}} 
\def\bs{\begin{small}}
\def\bS{\begin{slide}}
\def\ec{\end{center}} 
\def\ee{\end{enumerate}}
\def\ei{\end{itemize}}
\def\es{\end{small}}
\def\eS{\end{slide}}
\def\eeq{\end{eqnarray}}

\newcommand{\newP}[1]{\medskip\noindent{\bf #1:}}

\newcommand{\ud}{\,\mathrm{d}}

\def\Re{\mathbb{R}}
\def\E{{\sf E}}

\def\clY{{\cal Y}}

\def\clZ{{\cal Z}}

\renewcommand{\Re}{\mathbb{R}}